\documentclass[11pt,a4paper]{article}

\usepackage{inputenc}
\usepackage{amsmath}
\usepackage{bm}
\usepackage{bbold}
\usepackage{amsthm}

\usepackage{hyperref}

\title{Algebraic Solution to a Constrained Rectilinear Minimax Location Problem on the Plane\thanks{2011 International Conference on Multimedia Technology (ICMT), IEEE, 2011, pp.~6212--6220.}
} 

\author{Nikolai Krivulin\thanks{Faculty of Mathematics and Mechanics, St.~Petersburg State University, 28 Universitetsky Ave., St.~Petersburg, 198504, Russia, 
nkk@math.spbu.ru.} \thanks{The work was partially supported by the Russian Foundation for Basic Research under Grant \#09-01-00808.}
}

\date{}

\newtheorem{theorem}{Theorem}
\newtheorem{lemma}[theorem]{Lemma}
\newtheorem{corollary}[theorem]{Corollary}

\begin{document}

\maketitle

\begin{abstract}
We consider a constrained minimax single facility location problem on the plane with rectilinear distance. The feasible set of location points is restricted to rectangles with sides oriented at a 45 degrees angle to the axes of Cartesian coordinates. To solve the problem, an algebraic approach based on an extremal property of eigenvalues of irreducible matrices in idempotent algebra is applied. A new algebraic solution is given that reduces the problem to finding eigenvalues and eigenvectors of appropriately defined matrices.
\\

\end{abstract}

\section{Introduction}

Models and methods of idempotent algebra \cite{Cuninghame-Green1994Minimax,Kolokoltsov1997Idempotent,Golan2003Semirings,Heidergott2006Max-plus,Butkovic2010Maxlinear} find increasing application in solving actual problems in engineering, economics, management, and other fields. Expressed in terms of idempotent algebra, a range of problems that are nonlinear in the ordinary sense, become linear, resulting in more simple analysis and solution techniques. Specifically, many classical optimization problems (in particular, in graph optimization and dynamic programming) reduce to solving linear vector equations, finding eigenvalues and eigenvectors of matrices, and to working out similar computational procedures. 

Among the problems that have solutions in terms of idempotent algebra are certain location problems. A single facility one-dimensional location problem on a graph is examined in \cite{Cuninghame-Green1991Minimax,Cuninghame-Green1994Minimax}, where it is turned into a problem of minimizing a rational function in the idempotent algebra sense. Note, however, that the proposed solution concentrates on the analysis of polynomial and rational functions of one variable, and becomes less or no applicable in the case of two or more variables.

In \cite{Zimmermann2003Disjunctive,Tharwat2010Oneclass}, a multidimensional constrained location problem on a graph is reduced to an optimization problem with a max-separable objective function that can be represented as a maximum of functions each depending only on one variable. An efficient computational procedure is proposed which, however, seems to have limited application only to location problems with max-separable objective functions. 

In this paper, we consider a constrained rectilinear minimax single facility location problem on the plane, where the feasible set of location points is restricted to rectangles that have sides oriented at $45^{\circ}$ angle to the axes of Cartesian coordinates. In the unconstrained case, the problem has a well-known closed-form solution derived in \cite{Elzinga1972Geometrical,Francis1972Ageometrical} from geometrical arguments.

A new algebraic solution to the unconstrained problem is given in \cite{Krivulin2011Anextremal} on the basis of the approach proposed in \cite{Krivulin2011Algebraic,Krivulin2011Anextremal} to solve minimax location problems with rectilinear and Chebyshev metrics. The approach exploits properties of eigenvalues of irreducible matrices in idempotent algebra to reduce certain location problems to finding eigenvalues and eigenvectors of appropriately defined matrices. The aim of the paper is to extend the algebraic solution to the constrained location problem under consideration.

We start with a brief overview of basic facts about idempotent algebra to provide a background for subsequent sections. Some relevant results from the spectral theory of matrices are outlined including an extremal property of eigenvalues of irreducible matrices. We provide a conventional description for the location problem of interest and then represent the problem in terms of idempotent algebra. An algebraic solution is first given to the unconstrained problem and then extended to the constrained case.

\section{Preliminary Results}

In this section, a short introduction to idempotent algebra based on \cite{Krivulin2006Solution,Krivulin2006Eigenvalues} is given as background material to the solution approach. Further detail on the idempotent algebra theory and methods is provided in \cite{Cuninghame-Green1994Minimax,Kolokoltsov1997Idempotent,Golan2003Semirings,Heidergott2006Max-plus,Butkovic2010Maxlinear}.

\subsection{Idempotent Semifield}

Consider a set $\mathbb{X}$ endowed with two operations, addition $\oplus$ and multiplication $\otimes$. We assume $\langle\mathbb{X},\mathbb{0},\mathbb{1},\oplus,\otimes\rangle$ to be a commutative semiring with zero $\mathbb{0}$ and identity $\mathbb{1}$, where addition is idempotent and multiplication is invertible. Since the nonzero elements in $\mathbb{X}$ form a group with respect to multiplication, this semiring is often referred to as the idempotent semifield.

The power notation is defined as usual. Let $\mathbb{X}_{+}=\mathbb{X}\setminus\{\mathbb{0}\}$. For each $x\in\mathbb{X}_{+}$ and any integer $p>0$, we have $x^{0}=\mathbb{1}$, $\mathbb{0}^{p}=\mathbb{0}$, $x^{p}=x^{p-1}\otimes x=x\otimes x^{p-1}$, and $x^{-p}=(x^{-1})^{p}$.

We assume that the integer power can naturally be extended to the cases of rational and real exponents.

In the algebraic expressions below, the multiplication sign $\otimes$, as usual, is omitted. The power notation is thought in the sense of idempotent algebra.

The idempotent addition induces a partial order on $\mathbb{X}$ so that $x\leq y$ if and only if $x\oplus y=y$. With this definition, it is easy to verify that $x\leq x\oplus y$, $y\leq x\oplus y$, and that both addition and multiplication are isotonic.

Finally, it is not difficult to show that the binomial identity takes the form $(x\oplus y)^{\alpha}=x^{\alpha}\oplus y^{\alpha}$ for all $\alpha\geq0$.

As an example of such semirings, one can consider the idempotent semifield of real numbers 
$$
\mathbb{R}_{\max,+}
=
\langle\mathbb{R}\cup\{-\infty\},-\infty,0,\max,+\rangle.
$$

In $\mathbb{R}_{\max,+}$, the null element is $-\infty$, and the identity is $0$. For each $x\in\mathbb{R}$, there exists the inverse $x^{-1}$ that is equal to $-x$ in the conventional arithmetic. For any $x,y\in\mathbb{R}$, the notation $x^{y}$ conforms to the arithmetic product $xy$. The partial order coincides with the natural linear order in $\mathbb{R}$.

\subsection{Vector and Matrix Algebra}

Vector and matrix operations are routinely introduced on the basis of the scalar operations defined on $\mathbb{X}$. Consider a Cartesian product $\mathbb{X}^{n}$ with its elements represented as column vectors. For any vectors $\bm{x}=(x_{i})$ and $\bm{y}=(y_{i})$, and a scalar $c$, vector addition and multiplication by scalars follow the conventional rules
$$
\{\bm{x}\oplus\bm{y}\}_{i}
=
x_{i}\oplus y_{i},
\qquad
\{c\bm{x}\}_{i}
=
cx_{i}.
$$ 

%The set $\mathbb{X}^{n}$ equipped with these operations forms a vector semimodule over the idempotent semifield $\mathbb{X}$.

As usual, a vector $\bm{y}\in\mathbb{X}^{n}$ is linearly dependent on vectors $\bm{x}_{1},\ldots,\bm{x}_{m}\in\mathbb{X}^{n}$ if there are scalars $c_{1},\ldots,c_{m}\in\mathbb{X}$ such that $\bm{y}=c_{1}\bm{x}_{1}\oplus\cdots\oplus c_{m}\bm{x}_{m} $. In particular, $\bm{y}$ is collinear with $\bm{x}$ when $\bm{y}=c\bm{x}$.

For any column vector $\bm{x}=(x_{i})\in\mathbb{X}_{+}^{n}$, we define a row vector $\bm{x}^{-}=(x_{i}^{-})$ with elements $x_{i}^{-}=x_{i}^{-1}$. For all $\bm{x},\bm{y}\in\mathbb{X}_{+}^{n}$, the component-wise inequality $\bm{x}\leq\bm{y}$ leads to the inequality $\bm{x}^{-}\geq\bm{y}^{-}$.

For any conforming matrices $A=(a_{ij})$, $B=(b_{ij})$, and $C=(c_{ij})$ with entries in $\mathbb{X}$, matrix addition and multiplication together with multiplication by a scalar $c$ are performed in accordance with the formulas
$$
\{A\oplus B\}_{ij}
=
a_{ij}\oplus b_{ij},
\qquad
\{B C\}_{ij}
=
\bigoplus_{k}b_{ik}c_{kj},
\qquad
\{cA\}_{ij}=ca_{ij}.
$$

A matrix with all entries equal to zero is called the zero matrix and denoted by $\mathbb{0}$.

As is customary, a square matrix that has the diagonal entries equal to $\mathbb{1}$, and all off-diagonal entries equal to $\mathbb{0}$ is called the identity matrix and denoted by $I$. 

%The set $\mathbb{X}^{n\times n}$ with the matrix addition and multiplication is an idempotent semiring with identity.

For any square matrix $A\ne\mathbb{0}$ and integer $p>0$, we have $A^{0}=I $, $A^{p}=A^{p-1}A=AA^{p-1}$. The trace of a matrix $A=(a_{ij})\in\mathbb{X}^{n\times n}$ is given by $\mathop\mathrm{tr}A=a_{11}\oplus\cdots\oplus a_{nn}$.

A matrix is irreducible if it cannot be put into a block triangular form by simultaneous permutations of its rows and collumns.

\subsection{Eigenvalues and Eigenvectors}

A scalar $\lambda$ is an eigenvalue of a matrix $A\in\mathbb{X}^{n\times n}$ if there exists a vector $\bm{x}\in\mathbb{X}^{n}\setminus\{\mathbb{0}\}$ to provide the equality $A\bm{x}=\lambda\bm{x}$. Any vector $\bm{x}$ that satisfies this equality is an eigenvector of $A$ with the eigenvalue $\lambda$.

Any irreducible matrix $A$ has one eigenvalue given by
\begin{equation}
\lambda
=
\bigoplus_{k=1}^{n}\mathop\mathrm{tr}\nolimits^{1/k}(A^{k}).
\label{E-lambda}
\end{equation}

The eigenvectors of the matrix $A$ have no zero entries.

For any irreducible matrix $A\in\mathbb{X}^{n\times n}$, the eigenvectors corresponding to $\lambda$ can be found as follows. First we evaluate the matrix 
$$
A^{\times}
=
\lambda^{-1}A
\oplus
\cdots
\oplus
\lambda^{-n}A^{n}.
$$

Let $\bm{a}_{i}^{\times}$ be column $i$ in $A^{\times}$, and $a_{ii}^{\times}$ its diagonal element. We take the set of all columns $\bm{a}_{i}^{\times}$ such that $a_{ii}^{\times}=\mathbb{1}$, and then reduce it by removing those elements, if any, that are linearly dependent on others. Finally, the columns in the reduced set are put together to form a matrix $A^{+}$.  

The set of eigenvectors that corresponds to $\lambda$ (including the zero vector) coincides with the linear span of the columns of $A^{+}$, whereas any eigenvector is represented as  
$$
\bm{x}
=
A^{+}\bm{s},
$$
where $\bm{s}$ is a nonzero vector of appropriate size.

We conclude with the results from \cite{Krivulin2011Anextremal,Krivulin2011Algebraic} that underlies the solutions in the subsequent sections.
\begin{lemma}\label{L-mxmAx}
Let $A=(a_{ij})\in\mathbb{X}^{n\times n}$ be an irreducible matrix with an eigenvalue $\lambda$. Suppose $\bm{u}=(u_{i})$ and $\bm{v}=(v_{i})$ are eigenvectors of the matrices $A$ and $A^{T}$.

Then it holds that
$$
\min_{\bm{x}\in\mathbb{X}_{+}^{n}}\bm{x}^{-} A\bm{x}
=
\lambda,
$$
where the minimum is attained at any vector
\begin{equation}
\bm{x}
=
\left(
\begin{array}{c}
u_{1}^{\alpha}v_{1}^{\alpha-1} \\
\vdots \\
u_{n}^{\alpha}v_{n}^{\alpha-1}
\end{array}
\right),
\qquad
0\leq\alpha\leq1.
\label{E-xalpha}
\end{equation}
\end{lemma}

\section{A Minimax Location Problem}

In this section, we give two equivalent representations of the constrained minimax single facility location problem on a plane with rectilinear distance. We start with a problem definition in the usual notation, and then suggest an alternative form written in terms of the semifield $\mathbb{R}_{\max,+}$.

\subsection{Conventional Representation}

The distance between any two vectors $\bm{r}=(r_{1},r_{2})^{T}$ and $\bm{s}=(s_{1},s_{2})^{T}$ on the plane $\mathbb{R}^{2}$ in rectilinear metric is
\begin{equation}
\rho(\bm{r},\bm{s})
=
|r_{1}-s_{1}|+|r_{2}-s_{2}|.
\label{E-Manhattan}
\end{equation}

Suppose there are vectors $\bm{r}_{i}=(r_{1i},r_{2i})^{T}\in\mathbb{R}^{2}$ and numbers $w_{i}\in\mathbb{R}$ for $i=1,\ldots,m$, where $m\geq2$. Given a feasible set $S\subset\mathbb{R}^{2}$, the location problem of interest is to find the vectors $\bm{x}=(x_{1},x_{2})^{T}$ that provide for
\begin{equation}
\min_{\bm{x}\in S}\max_{1\leq i\leq m}(\rho(\bm{r}_{i},\bm{x})+w_{i}).
\label{P-Manhattan}
\end{equation}

In this paper, we restrict our consideration to feasible location sets that have the form of rectangles with sides oriented at $45^{\circ}$ angle to the axes of Cartesian coordinates.

%Note that feasible location regions of the above shape naturally arise when there are constraints on the maximum rectilinear distance allowed between each given point $\bm{r}_{i}$ and the facility location point $\bm{x}$. Since every constraint determines a rectangle turned by $45^{\circ}$, all constraints produce an intersection of the rectangles, which, if nonempty, has the same form.

\subsection{Algebraic Representation}

First we represent the metric \eqref{E-Manhattan} in terms of operations in the semifield $\mathbb{R}_{\max,+}$ as
$$
\rho(\bm{r},\bm{s})
=
(s_{1}^{-1}r_{1}\oplus r_{1}^{-1}s_{1})
(s_{2}^{-1}r_{2}\oplus r_{2}^{-1}s_{2}).
$$

Now problem \eqref{P-Manhattan} can be rewritten in the form
\begin{equation}
\min_{\bm{x}\in S}\varphi(\bm{x}),
\label{P-Manhattan1}
\end{equation}
where the objective function is given by
$$
\varphi(\bm{x})
=
\bigoplus_{i=1}^{m}w_{i}
(x_{1}^{-1}r_{1i}\oplus r_{1i}^{-1}x_{1})
(x_{2}^{-1}r_{2i}\oplus r_{2i}^{-1}x_{2}).
$$

To solve this constrained location problem we first examine its unconstrained version.

\section{The Unconstrained Problem}

Let us consider problem \eqref{P-Manhattan} and put $S=\mathbb{R}^{2}$ to get the problem
\begin{equation}
\min_{\bm{x}\in\mathbb{R}^{2}}\varphi(\bm{x}),
\label{P-Manhattan2}
\end{equation}

For this unconstrained problem, a solution based on geometrical arguments is obtained in \cite{Elzinga1972Geometrical}, where the problem is called the rectilinear messenger boy problem. A similar solution in the case when $w_{i}=0$ for all $i$ is given in \cite{Francis1972Ageometrical}.

Now we outline how the problem can be rearranged and then solved algebraically. Additional detail is given in \cite{Krivulin2011Anextremal}.

\subsection{Further Algebraic Transformation}

Consider the objective function of the problem. After some simple algebra, we get
\begin{equation}
\varphi(\bm{x})
=
ax_{1}^{-1}x_{2}
\oplus
bx_{1}x_{2}^{-1}
\oplus
cx_{1}^{-1}x_{2}^{-1}
\oplus
dx_{1}x_{2},
\label{E-phix}
\end{equation}
where
\begin{align}
a
&=
\bigoplus_{i=1}^{m}w_{i}r_{1i}r_{2i}^{-1},
&
b
&=
\bigoplus_{i=1}^{m}w_{i}r_{1i}^{-1}r_{2i},
\label{E-ab}
\\
c
&=
\bigoplus_{i=1}^{m}w_{i}r_{1i}r_{2i},
&
d
&=
\bigoplus_{i=1}^{m}w_{i}r_{1i}^{-1}r_{2i}^{-1}.
\label{E-cd}
\end{align}

Let us introduce a vector $\bm{y}$ and a matrix $A$ of third order as follows
\begin{equation}
\bm{y}
=
\left(
\begin{array}{c}
x_{1} \\
x_{2} \\
x_{1}^{-1}
\end{array}
\right),
\qquad
A
=
\left(
\begin{array}{ccc}
\mathbb{0} & a & \mathbb{0} \\
b & \mathbb{0} & c \\
\mathbb{0} & d & \mathbb{0}
\end{array}
\right).
\label{E-yA}
\end{equation}

Now the objective function can be represented as
$$
\varphi(\bm{x})
=
\bm{y}^{-}A\bm{y},
$$
whereas problem \eqref{P-Manhattan} reduces to that of the form
\begin{equation}
\min_{\bm{y}\in\mathbb{R}^{3}}\bm{y}^{-}A\bm{y}.
\label{P-minyAy}
\end{equation}

It is easy to see that the solutions of the new extended three-dimensional problem do not all correspond to solutions of \eqref{P-Manhattan2}. By virtue of the definitions at \eqref{E-yA}, a proper solution of the extended problem is a vector $\bm{y}$ such that its first and last components are inverse to each other.

\subsection{Algebraic Solution}

Let us examine extended problem \eqref{P-minyAy}. First note that the matrix $A$ defined in \eqref{E-yA} is irreducible. From Lemma~\ref{L-mxmAx} it follows that the minimum in \eqref{P-minyAy} is equal to the eigenvalue $\lambda$ of $A$, and the minimum is attained at any vector defined according to \eqref{E-xalpha} through eigenvectors of the matrices $A$ and $A^{T}$. Application of \eqref{E-lambda} gives
$$
\lambda
=
(ab\oplus cd)^{1/2}.
$$

Now we obtain the corresponding eigenvectors of the matrices $A$ and $A^{T}$. Note that $A^{T}$ is derived from $A$ by interchanging $a$ and $b$, as well as $c$ and $d$. Therefore, we can first find the eigenvector of $A$, and then turn it into the eigenvector of $A^{T}$ by the above interchange. 

To get the eigenvectors of $A$, we need to examine the matrix
$$
A^{\times}
=
\bigoplus_{k=1}^{3}\lambda^{-k}A^{k}
=
\left(
\begin{array}{ccc}
\lambda^{-2}ab & \lambda^{-1}a & \lambda^{-2}ac \\
\lambda^{-1}b & \mathbb{1} & \lambda^{-1}c \\
\lambda^{-2}bd & \lambda^{-1}d & \lambda^{-2}cd
\end{array}
\right).
$$

There are three cases to investigate.

\subsubsection{The case $\lambda^{2}=ab>cd$}

The matrix $A^{\times}$ takes the form
$$
A^{\times}
=
\left(
\begin{array}{ccc}
\mathbb{1} & \lambda^{-1}a & \lambda^{-2}ac \\
\lambda^{-1}b & \mathbb{1} & \lambda^{-1}c \\
\lambda^{-2}bd & \lambda^{-1}d & \lambda^{-2}cd
\end{array}
\right).
$$

Since $\lambda^{-2}cd<\mathbb{1}$, only the first two columns of $A^{\times}$ are eigenvectors of $A$. Considering that both columns are collinear, we take only one of them, say the first column, to form the matrix $A^{+}$. Now we arrive at the eigenvectors of $A$ defined as
$$
\bm{u}
=
A^{+}
s
=
\left(
\begin{array}{c}
\mathbb{1} \\
\lambda^{-1}b \\
\lambda^{-2}bd
\end{array}
\right)s,
\quad
s\in\mathbb{R}.
$$

By interchanging $a$ and $b$, and then $c$ and $d$, we get the eigenvectors of $A^{T}$
$$
\bm{v}
=
\left(
\begin{array}{c}
\mathbb{1} \\
\lambda^{-1}a \\
\lambda^{-2}ac
\end{array}
\right)t,
\quad
t\in\mathbb{R}.
$$

Application of \eqref{E-xalpha} leads to the solution of the extended problem in the form
$$
\bm{y}
=
\left(
\begin{array}{c}
\mathbb{1} \\
\lambda^{1-2\alpha}a^{\alpha-1}b^{\alpha} \\
\lambda^{2(1-2\alpha)}a^{\alpha-1}b^{\alpha}c^{\alpha-1}d^{\alpha}
\end{array}
\right)s^{\alpha}t^{\alpha-1},
\quad
s,t\in\mathbb{R}.
$$

Under the condition that the first and last components of the vector $\bm{y}$ must be reciprocal, we have 
$$
s^{\alpha}t^{\alpha-1}
=
(a^{\alpha}b^{\alpha-1}c^{1-\alpha}d^{-\alpha})^{1/2}.
$$

Turning back to problem \eqref{P-Manhattan}, we arrive at
$$
\bm{x}
=
\left(
\begin{array}{c}
(a^{\alpha}b^{\alpha-1}c^{1-\alpha}d^{-\alpha})^{1/2} \\
(a^{\alpha-1}b^{\alpha}c^{1-\alpha}d^{-\alpha})^{1/2}
\end{array}
\right).
$$

\subsubsection{The case $\lambda^{2}=cd>ab$}

In a similar way as above, we get the matrices
$$
\setlength{\arraycolsep}{4pt}
A^{\times}
=
\left(
\begin{array}{ccc}
\lambda^{-2}ab & \lambda^{-1}a & \lambda^{-2}ac \\
\lambda^{-1}b & \mathbb{1} & \lambda^{-1}c \\
\lambda^{-2}bd & \lambda^{-1}d & \mathbb{1}
\end{array}
\right)\!,
\quad
A^{+}
=
\left(
\begin{array}{c}
\lambda^{-2}ac \\
\lambda^{-1}c \\
\mathbb{1}
\end{array}
\right).
$$

The eigenvectors of $A$ and $A^{T}$ take the form
$$
\bm{u}
=
\left(
\begin{array}{c}
\lambda^{-2}ac \\
\lambda^{-1}c \\
\mathbb{1}
\end{array}
\right)s,
\quad
\bm{v}
=
\left(
\begin{array}{c}
\lambda^{-2}bd \\
\lambda^{-1}d \\
\mathbb{1}
\end{array}
\right)t,
\quad
s,t\in\mathbb{R}.
$$

The extended problem has the solution
$$
\bm{y}
=
\left(
\begin{array}{c}
\lambda^{2(1-2\alpha)}a^{\alpha}b^{\alpha-1}c^{\alpha}d^{\alpha-1} \\
\lambda^{1-2\alpha}c^{\alpha}d^{\alpha-1} \\
\mathbb{1}
\end{array}
\right)s^{\alpha}t^{\alpha-1},
\quad
s,t\in\mathbb{R}.
$$

To get a solution of \eqref{P-Manhattan2}, we put
$$
s^{\alpha}t^{\alpha-1}
=
(a^{-\alpha}b^{1-\alpha}c^{\alpha-1}d^{\alpha})^{1/2},
$$
and then get
$$
\bm{x}
=
\left(
\begin{array}{c}
(a^{\alpha}b^{\alpha-1}c^{1-\alpha}d^{-\alpha})^{1/2} \\
(a^{-\alpha}b^{1-\alpha}c^{\alpha}d^{\alpha-1})^{1/2}
\end{array}
\right).
$$

\subsubsection{The case $\lambda^{2}=ab=cd$}

From the matrix $A^{\times}$ we successively get
$$
A^{+}
=
\left(
\begin{array}{c}
\mathbb{1} \\
\lambda^{-1}b \\
\lambda^{-2}bd
\end{array}
\right),
\quad
\bm{x}
=
\left(
\begin{array}{c}
(a^{\alpha}b^{\alpha-1}c^{1-\alpha}d^{-\alpha})^{1/2} \\
(a^{\alpha-1}b^{\alpha}c^{1-\alpha}d^{-\alpha})^{1/2}
\end{array}
\right).
$$

\subsection{Summary of Solutions}

It is not difficult to verify that the above three solutions can be combined into one to represent the result as follows.
\begin{lemma}\label{L-Manhattan}
With the notation \eqref{E-ab}-\eqref{E-cd}, the minimum in \eqref{P-Manhattan2} is equal to
$$
\lambda=(ab\oplus cd)^{1/2}
$$
and attained at the vector
$$
\bm{x}
=
\left(
\begin{array}{c}
(a^{\alpha}b^{\alpha-1}c^{1-\alpha}d^{-\alpha})^{1/2} \\
\lambda^{2\alpha-1}(a^{-\alpha}b^{1-\alpha}c^{1-\alpha}d^{-\alpha})^{1/2}
\end{array}
\right)
$$
for all $\alpha$ such that $0\leq\alpha\leq1$.
\end{lemma}

In terms of ordinary arithmetic operations, we have the following statement.
\begin{corollary}
Suppose that
\begin{align*}
a
&=
\max_{1\leq i\leq m}(w_{i}+r_{1i}-r_{2i}),
&
b
&=
\max_{1\leq i\leq m}(w_{i}-r_{1i}+r_{2i}),
\\
c
&=
\max_{1\leq i\leq m}(w_{i}+r_{1i}+r_{2i}),
&
d
&=
\max_{1\leq i\leq m}(w_{i}-r_{1i}-r_{2i}).
\end{align*}

Then the minimum in \eqref{P-Manhattan} is given by
$$
\lambda=\max(a+b,c+d)/2,
$$
and it is attained at the vector
$$
\bm{x}
=
\left(
\begin{array}{c}
\frac{\alpha}{2}(a-d)-\frac{1-\alpha}{2}(b-c) \\
(2\alpha-1)\lambda-\frac{\alpha}{2}(a+d)+\frac{1-\alpha}{2}(b+c)
\end{array}
\right)
$$
for all $\alpha$ such that $0\leq\alpha\leq1$.
\end{corollary}

%As it easy to see, the above result is consistent with that in \cite{Elzinga1972Geometrical,Francis1972Ageometrical}.

\section{Constrained Problems}

The main idea that underlies the solution of the constrained problem under consideration is to put it in an unconstrained form, and then take advantage of results in the previous section. 

First we take the unconstrained problem at \eqref{P-Manhattan2}, and represent it in a slightly different way that is more convenient for further incorporation of constraints. With subscripts added to the above notation, we have
\begin{align*}
a_{0}
&=
\bigoplus_{i=1}^{m}w_{i}r_{1i}r_{2i}^{-1},
&
b_{0}
&=
\bigoplus_{i=1}^{m}w_{i}r_{1i}^{-1}r_{2i},
\\
c_{0}
&=
\bigoplus_{i=1}^{m}w_{i}r_{1i}r_{2i},
&
d_{0}
&=
\bigoplus_{i=1}^{m}w_{i}r_{1i}^{-1}r_{2i}^{-1},
\end{align*}
and
$$
\lambda_{0}
=
(a_{0}b_{0}\oplus c_{0}d_{0})^{1/2}.
$$

Furthermore, we define an objective function
$$
\varphi_{0}(\bm{x})
=
\lambda_{0}^{-1}(a_{0}x_{1}^{-1}x_{2}
\oplus
b_{0}x_{1}x_{2}^{-1}
\oplus
c_{0}x_{1}^{-1}x_{2}^{-1}
\oplus
d_{0}x_{1}x_{2}).
$$

Now consider the problem
$$
\min_{\bm{x}\in\mathbb{R}^{2}}\varphi_{0}(\bm{x}).
$$

It is easy to see that this problem is equivalent to \eqref{P-Manhattan2} in the sense that they both have the same solution set. Note, however, that the minimum in the last problem is $\mathbb{1}=0$, whereas it is equal to $\lambda=\lambda_{0}$ in \eqref{P-Manhattan2}.

\subsection{Location in 45$^{\circ}$ Rotated Rectangles}

Given numbers $a_{1},b_{1},c_{1},d_{1}$, we assume the feasible location set $S$ to be the intersection of the half-planes defined in standard notation by the inequalities
\begin{align*}
-x_{1}+x_{2}+a_{1}
&\leq
0,
&
x_{1}-x_{2}+b_{1}
&\leq
0,
%\label{E-a1b1}
\\
-x_{1}-x_{2}+c_{1}
&\leq
0,
&
x_{1}+x_{2}+d_{1}
&\leq
0.
%\label{E-c1d1}
\end{align*}

It is easy to see that the intersection of the half-planes, if nonempty, has the form of an upright rectangle rotated by 45$^\circ$ around its center.

In terms of the semifield $\mathbb{R}_{\max,+}$, the inequalities can be rewritten as
\begin{align*}
a_{1}x_{1}^{-1}x_{2}
&\leq
\mathbb{1},
&
b_{1}x_{1}x_{2}^{-1}
&\leq
\mathbb{1},
%\label{E-a1b1}
\\
c_{1}x_{1}^{-1}x_{2}^{-1}
&\leq
\mathbb{1},
&
d_{1}x_{1}x_{2}
&\leq
\mathbb{1},
%\label{E-c1d1}
\end{align*}
and then further combined into one inequality
$$
a_{1}x_{1}^{-1}x_{2}
\oplus
b_{1}x_{1}x_{2}^{-1}
\oplus
c_{1}x_{1}^{-1}x_{2}^{-1}
\oplus
d_{1}x_{1}x_{2}
\leq
\mathbb{1}.
$$

Let us introduce a function
$$
\varphi_{1}(\bm{x})
=
a_{1}x_{1}^{-1}x_{2}
\oplus
b_{1}x_{1}x_{2}^{-1}
\oplus
c_{1}x_{1}^{-1}x_{2}^{-1}
\oplus
d_{1}x_{1}x_{2},
$$
and note that $\varphi_{1}(\bm{x})\leq\mathbb{1}$ if and only if $\bm{x}\in S$.

We define an objective function
$$
\psi(\bm{x})=\varphi_{0}(\bm{x})\oplus\varphi_{1}(\bm{x}).
$$

Furthermore, with the notation
\begin{align*}
a
&=
\lambda_{0}^{-1}a_{0}\oplus a_{1},
&
b
&=
\lambda_{0}^{-1}b_{0}\oplus b_{1},
\\
c
&=
\lambda_{0}^{-1}c_{0}\oplus c_{1},
&
d
&=
\lambda_{0}^{-1}d_{0}\oplus d_{1},
\end{align*}
we can put the objective function in the form of \eqref{E-phix}
$$
\psi(\bm{x})
=
ax_{1}^{-1}x_{2}
\oplus
bx_{1}x_{2}^{-1}
\oplus
cx_{1}^{-1}x_{2}^{-1}
\oplus
dx_{1}x_{2}.
$$

Now the constrained problem \eqref{P-Manhattan1} reduces to an unconstrained problem
\begin{equation}
\min_{\bm{x}\in\mathbb{R}^{2}}\psi(\bm{x})
\label{P-Manhattan3}
\end{equation}
that has the form of \eqref{P-Manhattan2} and the solution given by Lemma~\ref{L-Manhattan}.

\subsection{Solution to the Constrained Problem}

First note that the minimum in problem \eqref{P-Manhattan3} is given by $\lambda=(ab\oplus cd)^{1/2}\geq\mathbb{1}$. It is not difficult to see that $\lambda=\mathbb{1}$ if and only if the feasible set of the original constrained problem is nonempty and it has an intersection with the solution of the associated unconstrained problem.

In the case of $\lambda>\mathbb{1}$, the solution to \eqref{P-Manhattan3} can be considered as an approximate solution for the constrained problem, that is of independent interest.

In the conventional notation, the solution to \eqref{P-Manhattan3} can be described as follows.
\begin{corollary}
Suppose that
\begin{align*}
a
&=
\max(a_{0}-\lambda_{0},a_{1}),
&b
&=
\max(b_{0}-\lambda_{0},b_{1}),
\\
c
&=
\max(c_{0}-\lambda_{0},c_{1}),
&d
&=
\max(d_{0}-\lambda_{0},d_{1}),
\end{align*}
where
\begin{align*}
a_{0}
&=
\max_{1\leq i\leq m}(w_{i}+r_{1i}-r_{2i}),
\\
b_{0}
&=
\max_{1\leq i\leq m}(w_{i}-r_{1i}+r_{2i}),
\\
c_{0}
&=
\max_{1\leq i\leq m}(w_{i}+r_{1i}+r_{2i}),
\\
d_{0}
&=
\max_{1\leq i\leq m}(w_{i}-r_{1i}-r_{2i}),
\\
\lambda_{0}
&=
\max(a_{0}+b_{0},c_{0}+d_{0})/2.
\end{align*}

Then the minimum in \eqref{P-Manhattan3} is given by
$$
\lambda=\max(a+b,c+d)/2,
$$
and it is attained at the vector
$$
\bm{x}
=
\left(
\begin{array}{c}
\frac{\alpha}{2}(a-d)-\frac{1-\alpha}{2}(b-c) \\
(2\alpha-1)\lambda-\frac{\alpha}{2}(a+d)+\frac{1-\alpha}{2}(b+c)
\end{array}
\right)
$$
for all $\alpha$ such that $0\leq\alpha\leq1$.
\end{corollary}

\section*{Acknowledgment}
The work was partially supported by the Russian Foundation for Basic Research (grant No.~09-01-00808).

\bibliographystyle{utphys}

\bibliography{Algebraic_solution_to_a_constrained_rectilinear_minimax_location_problem_on_the_plane}

\end{document}